\documentclass[11pt,reqno]{amsart}
\usepackage{amsmath,amsfonts, amssymb, amsthm}
  \usepackage{paralist}
  \usepackage{graphics} 
  \usepackage{epsfig} 
\usepackage{graphicx}  \usepackage{epstopdf}
 \usepackage[colorlinks=true]{hyperref}
\hypersetup{urlcolor=blue, citecolor=red}

\newtheorem{theorem}{Theorem}[section]
\newtheorem{corollary}{Corollary}[section]

\theoremstyle{definition}

\newtheorem{remark}{Remark}[section]

\def\Z{\Bbb Z}
\def\N{\Bbb N}

\def\C{\Bbb C}

\def\l{\left}
\def\r{\right}
\def\bg{\bigg}
\def\({\bg(}
\def\){\bg)}
\def\t{\text}
\def\f{\frac}

\def\ls{\leq}
\def\gs{\geq}

\def\sm{\setminus}

\def\bs{\boldsymbol}
\def\al{\alpha}

\begin{document}
\hbox{Preprint, {\tt arXiv:2303.12629}}
\medskip

\title[A new result similar to the Graham-Pollak theorem]
      {A new result similar to\\ the Graham-Pollak theorem}
\author[Zhi-Wei Sun]{Zhi-Wei Sun}


\address{Department of Mathematics, Nanjing
University, Nanjing 210093, People's Republic of China}
\email{{\tt zwsun@nju.edu.cn}
\newline\indent
{\it Homepage}: {\tt http://maths.nju.edu.cn/\lower0.5ex\hbox{\~{}}zwsun}}

\keywords{Determinant, distance, $q$-analogue, tree.
\newline \indent Supported by the Natural Science Foundation of China (grant no. 11971222).}

\begin{abstract}
Let $n>1$ be an integer, and let $T$ be a tree with $n+1$ vertices $v_1,\ldots,v_{n+1}$, where
$v_1$ and $v_{n+1}$ are two leaves of $T$. For each edge $e$ of $T$, assign a complex number $w(e)$ as its weight. We obtain that
$$\det[x+d(v_{j+1},v_k)]_{1\ls j,k\ls n}=2^{n-2}\prod_{e\in E(T)}w(e),$$
where $d(v_{j+1},v_k)$ is the weighted distance between $v_{j+1}$ and $v_k$ in the tree $T$.
This is similar to the celebrated Graham-Pollak theorem on determinants of distance matrices for trees.
Actually, a more general result is deduced in this paper.
\end{abstract}
\maketitle

\section{Introduction}
\setcounter{equation}{0}
 \setcounter{conjecture}{0}
 \setcounter{theorem}{0}
 \setcounter{proposition}{0}

In 1971, R.L. Graham and H.O. Pollak \cite{GP} obtained the following beautiful result for (undirected) trees.

\begin{theorem}[Graham-Pollak Theorem] Let $T$ be any tree with $n$ vertices $v_1,\ldots,v_n$.
Then the determinant of the distance matrix $[d(v_j,v_k)]_{1\ls j,k\ls n}$
has the value $(1-n)(-2)^{n-2}$, where
$d(v_j,v_k)$ is the distance between $v_j$ and $v_k$ in the tree $T$.
\end{theorem}

For simple induction proofs of the Graham-Pollak theorem, one may consult Z. Du and J. Yeh \cite {DY}, and also W. Yan and Y.-N. Yeh \cite{Y06}. Note also that R. Robinson and G. Szeg\"o \cite{RS} evaluated $\det[|j-k|]_{1\ls j,k\ls n}$ in 1936.

The $q$-analogue of an integer $m$ is given by
$$[m]_q=\f{q^m-1}{q-1},$$
where we view $q$ as a variable.
Obviously $\lim_{q\to1}[m]_q=m$ for all $m\in\N$, and it is easy to verify that
\begin{equation}\label{m}[m_1+m_2]_q=q^{m_1}[m_2]_q+[m_1]_q\ \ \ \ \t{for all}\ m_1,m_2\in\Z.
\end{equation}
R.B. Bapat, A.K. Lal and S. Pati \cite{BLP}, as well as  Yan and Yeh \cite{YY}, obtained the $q$-analogue of the Grahm-Pollak theorem which states that for any tree with $n>1$ vertices $v_1,\ldots,v_n$ we have
\begin{equation}\label{q-tree}\det[[d(v_j,v_k)]_q]_{1\ls j,k\ls n}=(1-n)(-1-q)^{n-2}.
\end{equation}

Let $T$ be a weighted tree with $n>1$ vertices $v_1,\ldots,v_n$, and suppose that each edge $e$ of $T$ is assigned a weight $w(e)\in\Z$. In 2005, Bapat, S.J. Kirkland and M. Neuman \cite{BKN} deduced that
\begin{equation}\label{GPx}\det[x+d(v_j,v_k)]_{1\ls j,k\ls n}=(-1)^{n-1}2^{n-2}\(2x+\sum_{e\in E(T)}w(e)\)\prod_{e\in E(T)}w(e),
\end{equation}
where $E(T)$ is the edge set of $T$, and $d(v_j,v_k)$ is the weighted distance which is the sum of the weights of all edges in the unique path from $v_j$ to $v_k$.
Bapat and P. Rekhi \cite{BR} proved that
\begin{equation}\label{w-tree}\det[[d(v_j,v_k)]_q]_{1\ls j,k\ls n}=(-1)^{n-1}\(\prod_{e\in E(T)}[2w(e)]_q\)\sum_{e\in E(T)}\f{[w(e)]_q}{1+q^{w(e)}},
\end{equation}
i.e.,
\begin{equation}\label{w-t}\det[q^{d(v_j,v_k)}-1]_{1\ls j,k\ls n}=\(\prod_{e\in E(T)}(1-q^{2w(e)})\)\sum_{e\in E(T)}\f{q^{w(e)}-1}{q^{w(e)}+1}.
\end{equation}
 Yan and Yeh \cite{YY} got another formula for $\det[[d(v_j,v_k)]_q]_{1\ls j,k\ls n}$
involving weighted distances, which is more complicated than \eqref{w-tree}.
Bapat, Lal and Pati \cite{BLP} obtained that
\begin{equation}\label{q0}\det[q^{d(v_j,v_k)}]_{1\ls j,k\ls n}=\prod_{e\in E(T)}(1-q^{2w(e)}),\end{equation}
which was also deduced by Yan and Yeh \cite{YY}. As $\det[q^{d(v_j,v_k)}-t]_{1\ls j,k\ls n}$
is linear in $t$ (which can be easily seen), combining \eqref{w-t} and \eqref{q0} we find that
\begin{equation}\label{GPt}\det[q^{d(v_j,v_k)}-t]_{1\ls j,k\ls n}=\(\prod_{e\in E(T)}(1-q^{2w(e)})\)\(1+t\(\sum_{e\in E(T)}\f{q^{w(e)}-1}{q^{w(e)}+1}-1\)\).
\end{equation}

A pendant vertex (i.e., a vertex of degree one) in a tree is called a {\it leaf} of the tree. It is well known that any tree with $n>1$ vertices have at least two leaves and exactly $n-1$ edges.

In this paper we establish the following new theorem.

\begin{theorem}\label{Th1.1} Let $n>1$ be an integer, and let $T$ be a tree with $n+1$ vertices $v_1,\ldots,v_{n+1}$
and $n$ edges $e_1,\ldots,e_n$, where
$v_1$ and $v_{n+1}$ are two leaves of $T$, and $e_1$ and $e_n$ are incident to the leaves $v_1$ and $v_{n+1}$
respectively. For each edge $e$ of $T$, assign a weight $w(e)\in\Z$.
 Then
 \begin{equation} \label{qt}\det[q^{d(v_{j+1},v_k)}-t]_{1\ls j,k\ls n}=t(q^{w(e_1)}-1)(q^{w(e_n)}-1)
\prod_{1<i<n}(q^{2w(e_i)}-1).
\end{equation}
In particular,
\begin{equation} \label{wq} \det[x+d_q(v_{j+1},v_k)]_{1\ls j,k\ls n}=((1-q)x+1)[w(e_1)]_q[w(e_n)]_q\prod_{1<i<n}[2w(e_i)]_q,
\end{equation}
where $d_q(v_{j+1},v_k)$ denotes the $q$-analogue of the weighted distance $d(v_{j+1},v_k)$.
\end{theorem}
\begin{remark} Note that \eqref{wq} reduces to the identity
\begin{equation}\label{nq}\det[x+d_q(v_{j+1},v_k)]_{1\ls j,k\ls n}=((1-q)x+1)(1+q)^{n-2}
\end{equation}
provided that all the edges $e_1,\ldots,e_{n}$ have the usual weight $1$.
\end{remark}

For the identity \eqref{nq} with $x=0$, we still need the condition that $v_1$ and $v_{n+1}$ are two leaves of $T$.
This is illustrated by the following example.

{\it Example} 1.1. Consider a (not weighted) tree $T$ with vertices $v_1,v_2,v_3,v_4$ and edges $v_1v_2,v_1v_3,v_1v_4$.
Clearly $v_4$ is a leaf of $T$ but $v_1$ is not. Note that
$$\det[d_q(v_{j+1},v_k)]_{1\ls j,k\ls 3}=\vmatrix 1&0&[2]_q\\1&[2]_q&0\\1&[2]_q&[2]_q\endvmatrix=2[2]_q^2-[2]_q^2=(1+q)^2$$
is different from $(1+q)^{3-2}$.

Now we give an example to illustrate \eqref{nq} with $x=0$ and $n=7$.

{\it Example} 1.2. Let us consider a tree $T$ with $8$ vertices $v_1,\ldots,v_8$
and edges
$$v_1v_2,\ v_2v_3,\ v_3v_4,\ v_5v_6,\ v_6v_7,\ v_2v_7,\ v_7v_8.$$
Then $v_1$ and $v_8$ are two leaves. Note that the vertex $v_2$ adjacent to $v_1$
is neither of degree two nor adjacent to another leaf. Similarly, the vertex $v_7$ adjacent to $v_8$
is neither of degree two nor adjacent to another leaf. Observe that
\begin{align*}
&\ \det[d_q(v_{j+1},v_k)]_{1\ls j,k\ls 7}
\\=&\ \vmatrix 1&0&1&[2]_q&[3]_q&[2]_q&1\\ [2]_q&1&0&1&[4]_q&[3]_q&[2]_q
\\ [3]_q&[2]_q&1&0&[5]_q&[4]_q&[3]_q
\\ [4]_q&[3]_q&[4]_q&[5]_q&0&1&[2]_q
\\ [3]_q&[2]_q&[3]_q&[4]_q&1&0&1
\\ [2]_q&1&[2]_q&[3]_q&[2]_q&1&0
\\ [3]_q&[2]_q&[3]_q&[4]_q&[3]_q&[2]_q&1
\endvmatrix=(1+q)^{7-2}.
\end{align*}

Applying \eqref{wq} with $q\to1$, we get the following result similar to the Graham-Pollak theorem.

\begin{corollary}\label{Cor-w} Let $n>1$ be an integer, and let $T$ be a tree with $n+1$ vertices $v_1,\ldots,v_{n+1}$, where
$v_1$ and $v_{n+1}$ are two leaves of $T$. For each edge $e$ of $T$, assign a weight $w(e)\in\Z$. Then we have the identity
\begin{equation}\label{new}\det[x+d(v_{j+1},v_k)]_{1\ls j,k\ls n}=2^{n-2}\prod_{e\in E(T)}w(e).
\end{equation}
\end{corollary}
\begin{remark} Actually $w(e)\in\Z$ in Corollary \ref{Cor-w} can be replaced by $w(e)\in\C$, where
 $\C$ is the field of complex numbers; this is because any $\al\in\C$
 has its $q$-analogue $[\al]_q=(q^{\al}-1)/(q-1)$ provided that $q>0$ and $q\not=1$.
 Instead of using integer weights in Corollary \ref{Cor-w}, we may even utilize weights in any commutative ring with identity
; for this purpose we may modify our proof of Theorem \ref{Th1.1} by using weights instead of their $q$-analogues.
\end{remark}

If we apply Theorem \ref{Th1.1} to a path with vertices $v_1,\ldots,v_{n+1}$ and edges $v_iv_{i+1}\ (1\ls i\ls n)$, then we obtain the following result.

\begin{corollary} For any integer $n\gs2$, we have
\begin{equation}\det[q^{|j-k+1|}-t]_{1\ls j,k\ls n}=t(q-1)^n(q+1)^{n-2}.
\end{equation}
\end{corollary}

We are going to reduce Theorem \ref{Th1.1} to a new version in the next section, and then prove the new version in Section 3.

\section{A new version of Theorem \ref{Th1.1}}
 \setcounter{equation}{0}
 \setcounter{conjecture}{0}
 \setcounter{theorem}{0}
 \setcounter{proposition}{0}

Under the condition of Theorem \ref{Th1.1}, the identity \eqref{qt} with $t=(1-q)x+1$ yields \eqref{wq}.
Since
$$q^{d(v_{j+1},v_k)}-t-\l(q^{d(v_2,v_k)}-t\r)=q^{d(v_{j+1},v_k)}-q^{d(v_2,v_k)}$$
for all $1<j\ls n$ and $1\ls k\ls n$,
 we see that
$$\det[q^{d(v_{j+1},v_k)}-t]_{1\ls j,k\ls n}=ct+\det[q^{d(v_{j+1},v_k)}]_{1\ls j,k\ls n}$$
for certain $c$ independent of $t$. Let $v_{j+1}$ with $1\ls j\ls n-1$ be the unique vertex adjacent to the leaf $v_{n+1}$.
Then
$$q^{d(v_{n+1},v_k)}=q^{d(v_{j+1},v_k)+w(e_n)}=q^{w(e_n)}\times q^{d(v_{j+1},v_k)}$$
for all $k=1,\ldots,n$, and hence $\det[q^{d(v_{j+1},v_k)}]_{1\ls j,k\ls n}=0$.
Therefore,
$$\det[q^{d(v_{j+1},v_k)}-t]_{1\ls j,k\ls n}=ct$$
and $c=\det[q^{d(v_{j+1},v_k)}-1]_{1\ls j,k\ls n}$.
It follows that \eqref{qt} is equivalent to \eqref{wq} with $x=0$.
Note that $\det[d_q(v_{j+1},v_k)]_{1\ls j,k\ls n}$ does not depend on the order of $v_2,\ldots,v_n$.
In fact, if $1\ls s<t\ls n-1$ then by exchanging the $s$th row with the $t$th row, and also exchanging the $(s+1)$-th column with the $(t+1)$-th column, we obtain from $[d_q(v_{j+1},v_k)]_{1\ls j,k\ls n}$  the matrix $[d_q(u_{j+1},u_k)]_{1\ls j,k\ls n}$, where $u_{s+1}=v_{t+1}$, $u_{t+1}=v_{s+1}$ and $u_{i}=v_i$ for all $1\ls i\ls n+1$ with $i\not=s+1,t+1$.
So, without any loss of generality, we may assume that $v_2$ is the unique vertex adjacent to $v_1$, and also $v_n$ is the unique vertex adjacent to $v_{n+1}$. In view of \eqref{m},
\begin{align*}&d_q(v_{j+1},v_1)-q^{w(e_1)}d_q(v_{j+1},v_2)
\\=\ &[d(v_{j+1},v_2)+w(e_1)]_q-q^{w(e_1)}[d(v_{j+1},v_2)]_q=[w(e_1)]_q
\end{align*}
for all $j=1,\ldots,n$, and also
\begin{align*}&d_q(v_{n+1},v_k)-q^{w(e_n)}d_q(v_n,v_k)
\\=\ &[d(v_n,v_k)+w(e_n)]_q-q^{w(e_n)}[d(v_n,v_k)]_q=[w(e_n)]_q
\end{align*}
for all $k=1,\ldots,n$. So we have
\begin{align*}&\det[d_q(v_{j+1},v_k)]_{1\ls j,k\ls n}
\\=&\ \vmatrix [w(e_1)]_q&d_q(v_2,v_2)&d_q(v_2,v_3)&\hdots&d_q(v_2,v_n)\\ [w(e_1)]_q&d_q(v_3,v_2)&d_q(v_3,v_3)&\hdots&d_q(v_3,v_n)
\\\vdots&\vdots&\vdots&\ddots&\vdots\\ [w(e_1)]_q&d_q(v_n,v_2)&d_q(v_n,v_3)&\hdots&d_q(v_n,v_n)
\\(1-q^{w(e_n)})[w(e_1)]_q&[w(e_n)]_q&[w(e_n)]_q&\hdots&[w(e_n)]_q\endvmatrix
\\=&\ [w(e_1)]_q[w(e_n)]_q\vmatrix 1&d_q(v_2,v_2)&d_q(v_2,v_3)&\hdots&d_q(v_2,v_n)\\ 1&d_q(v_3,v_2)&d_q(v_3,v_3)&\hdots&d_q(v_3,v_n)
\\\vdots&\vdots&\vdots&\ddots&\vdots\\ 1&d_q(v_n,v_2)&d_q(v_n,v_3)&\hdots&d_q(v_n,v_n)
\\ 1-q&1&1&\hdots&1\endvmatrix.
\end{align*}
Note that, for the tree $T-\{v_1,v_{n+1}\}$, the vertices are $v_2,\ldots,v_n$
and the edges are those $e_i$ with $1<i<n$.
Therefore, Theorem \ref{Th1.1} has the following equivalent version.

\begin{theorem}\label{Th2.1} Let $T$ be any tree with vertices $v_1,\ldots,v_n$ and the edge set $E(T)$.
For each $e\in E(T)$ assign a weight $w(e)\in\Z$. Define
\begin{equation}\label{Dn} D_T(q):=\vmatrix 1&d_q(v_1,v_1)&d_q(v_1,v_2)&\hdots&d_q(v_1,v_n)\\1&d_q(v_2,v_1)&d_q(v_2,v_2)&\hdots&d_q(v_2,v_n)
\\\vdots&\vdots&\vdots&\ddots&\vdots\\1&d_q(v_n,v_1)&d_q(v_n,v_2)&\hdots&d_q(v_n,v_n)
\\1-q&1&1&\hdots&1\endvmatrix,\end{equation}
where $d_q(v_j,v_k)$ is the $q$-analogue of the weighted distance $d(v_j,v_k)$.
Then we have
\begin{equation}D_T(q)=\prod_{e\in E(T)}[2w(e)]_q.\end{equation}
\end{theorem}

\section{Proof of Theorem \ref{Th2.1}}
 \setcounter{equation}{0}
 \setcounter{conjecture}{0}
 \setcounter{theorem}{0}
 \setcounter{proposition}{0}

\medskip
 \noindent{\it Proof of Theorem \ref{Th2.1}}.
 Note that $D_T(q)$ keeps unchanged if we relabel the vertices of the tree $T$.
Below we prove the desired result by induction on $n=|V(T)|$, where $V(T)$ is the vertex set of $T$.

 If $n=1$, then $E(T)=\emptyset$, and
 $$D_T(q)=\vmatrix 1&0\\1-q&1
\endvmatrix=1=\prod_{e\in E(T)}[2w(e)]_q$$
When $n=2$, clearly $e=v_1v_2$ is the unique edge of $T$, and
\begin{align*}D_T(q)=&\ \vmatrix 1&0&[w(e)]_q\\1&[w(e)]_q&0\\1-q&1&1\endvmatrix
\\=\ &2[w(e)]_q-(1-q)[w(e)]_q^2=(q^{w(e)}+1)[w(e)]_q=[2w(e)]_q.
\end{align*}
So Theorem \ref{Th2.1} holds for $n=1,2$

Now we consider the case $n=3$. In this case, the tree $T$ is a path, and we may assume that $e_1=v_1v_2$ and $e_2=v_2v_3$ are its edges. Thus
\begin{align*}D_T(q)=&\ \vmatrix 1&d_q(v_1,v_1)&d_q(v_1,v_2)&d_q(v_1,v_3)\\1&d_q(v_2,v_1)&d_q(v_2,v_2)&d_q(v_2,v_3)
\\1&d_q(v_3,v_1)&d_q(v_3,v_2)&d_q(v_3,v_3)\\1-q&1&1&1\endvmatrix
\\=\ & \vmatrix 1&0&[w(e_1)]_q&[w(e_1)+w(e_2)]_q\\1&[w(e_1)]_q&0&[w(e_2)]_q
\\1&[w(e_1)+w(e_2)]_q&[w(e_2)]_q&0\\1-q&1&1&1\endvmatrix
\end{align*}
It is easy to verify that this determinant has the value $[2w(e_1)]_q[2w(e_2)]_q$ as desired.

Now, let $n=|V(T)|\gs4$ and assume the desired result for trees with smaller number of vertices. For any vector ${\bs m}=(m_1,\ldots,m_{n-3})$ with $m_1,\ldots,m_{n-3}\in\Z$, we simply write the vector
$([m_1]_q,\ldots,[m_{n-3}]_q)$ as $[{\bs m}]_q$, and let $[{\bs m}]_q'$ denote the column vector
which is the transpose of $[{\bs m}]_q$. For convenience, we also use $\bf 0$ and ${\bf 1}$ to denote the vectors
$(0,\ldots,0)$ and $(1,\ldots,1)$ of length $n-3$, respectively.

Without loss of generality, we may assume that $v_n$ is a terminal vertex of a path with the greatest length.
Note that the unique vertex adjacent to $v_n$ is of degree two or adjacent to another leaf.
Thus, without loss of generality, we may also suppose that $v_{n-1}$
is a vertex only adjacent to $v_n$ and $v_{n-2}$, or another leaf sharing a common neighbor $v_{n-2}$ with $v_n$.
Let $e_1$ stand for the edge $v_{n-1}v_{n-2}$, and let ${\bs u}$ denote the row vector
$$(d(v_{n-2},v_1),\ldots,d(v_{n-2},v_{n-3}))=(d(v_1,v_{n-2}),\ldots,d(v_{n-3},v_{n-2})).$$
It is easy to see that
\begin{equation}\label{n-2}\bmatrix [d_q(v_j,v_k)]_{1\ls j,k\ls n-3}&[{\bs u}]_q'\\ [{\bs u}]_q&0
\endbmatrix=[d_q(v_j,v_k)]_{1\ls j,k\ls n-2}
\end{equation}
and
\begin{equation}\label{n-1}\bmatrix [d_q(v_j,v_k)]_{1\ls j,k\ls n-3}&[{\bs u}]_q'&[{\bs u}+a{\bf 1}]_q'\\ [{\bs u}]_q&0&[a]_q
\\ [{\bs u}+a{\bf 1}]_q&[a]_q&0
\endbmatrix=[d_q(v_j,v_k)]_{1\ls j,k\ls n-1},
\end{equation}
where $a$ stands for $w(e_1)$.
\medskip

{\it Case} 1. $v_{n-1}$ is only adjacent to $v_n$ and $v_{n-2}$.
\medskip

Let $e_2$ denote the edge $v_{n-1}v_n$, and write $b$ for $w(e_2)$. Then
$$D_T(q)=\vmatrix {\bf 1}'&[d_q(v_j,v_k)]_{1\ls j,k\ls n-3}&[{\bs u}]_q'&[{\bs u}+a{\bf 1}]_q'&[{\boldsymbol u}+
(a+b){\bf 1}]_q'\\1&[{\bs u}]_q&0&[a]_q&[a+b]_q
\\1&[{\bs u}+a{\bf 1}]_q&[a]_q&0&[b]_q
\\1&[{\bs u}+(a+b){\bf 1}]_q&[a+b]_q&[b]_q&0
\\1-q&{\bf 1}&1&1&1\endvmatrix.$$
Replacing the $n$th row $R_n$ in $D_T(q)$ by $R_n-q^b\times R_{n-1}-[b]_q\times R_{n+1}$, we get that
$$D_T(q)=\vmatrix {\bf 1}'&[d_q(v_j,v_k)]_{1\ls j,k\ls n-3}&[{\bs u}]_q'&[{\bs u}+a{\bf 1}]_q'&[{\boldsymbol u}+(a+b){\bf 1}]_q'\\1&[{\bs u}]_q&0&[a]_q&[a+b]_q
\\1&[{\bs u}+a{\bf 1}]_q&[a]_q&0&[b]_q
\\0&{\bf 0}&0&0&-(q^b+1)[b]_q
\\1-q&{\bf 1}&1&1&1\endvmatrix$$
with the aid of \eqref{m}.
 Expanding this determinant according to its $n$th row, we find that
 \begin{align*}D_T(q)=&\ -(-(q^b+1)[b]_q)
 \vmatrix {\bf 1}'&[d_q(v_j,v_k)]_{1\ls j,k\ls n-3}&[{\bs u}]_q'&[{\bs u}+a{\bf 1}]_q'\\1&[{\bs u}]_q&0&[a]_q
\\1&[{\bs u}+a{\bf 1}]_q&[a]_q&0
\\1-q&{\bf 1}&1&1\endvmatrix
\\=\ &[2b]_qD_{T-\{v_n\}}(q)
\end{align*} in view of \eqref{n-1}. Thus, by the induction hypothesis for the tree $T-\{v_n\}$, we have
$$D_T(q)=[2b]_q\times \prod_{e\in E(T-\{v_n\})}[2w(e)]_q=\prod_{e\in E(T)}[2w(e)]_q.$$

{\it Case} 2. $v_{n-1}$ is another leaf sharing a common neighbor $v_{n-2}$ with $v_n$.
\medskip

In this case, we let $e_2$ be the edge $v_{n-2}v_n$, and write $b$ for $w(e_2)$. Then
$$D_T(q)=\vmatrix {\bf 1}'&[d_q(v_j,v_k)]_{1\ls j,k\ls n-3}&[{\bs u}]_q'&[{\bs u}+a{\bf 1}]_q'&[{\boldsymbol u}+b{\bf 1}]_q'\\1&[{\bs u}]_q&0&[a]_q&[b]_q
\\1&[{\bs u}+a{\bf 1}]_q&[a]_q&0&[a+b]_q
\\1&[{\bs u}+b{\bf 1}]_q&[b]_q&[a+b]_q&0
\\1-q&{\bf 1}&1&1&1\endvmatrix.$$
Replacing the $n$th row $R_n$ by $R_n-q^{b-a}\times R_{n-1}-[b-a]_q\times R_{n+1}$, we find that $D_T(q)$
equals
$$\vmatrix {\bf 1}'&[d_q(v_j,v_k)]_{1\ls j,k\ls n-3}&[{\bs u}]_q'&[{\bs u}+a{\bf 1}]_q'&[{\boldsymbol u}+b{\bf 1}]_q'\\1&[{\bs u}]_q&0&[a]_q&[b]_q
\\1&[{\bs u}+a{\bf 1}]_q&[a]_q&0&[a+b]_q
\\0&{\bf 0}&0&f(a,b)&f(a,b)-(1+q^{b-a})[a+b]_q
\\1-q&{\bf 1}&1&1&1\endvmatrix$$
with the aid of \eqref{m}, where
$$f(a,b)=[a+b]_q-[b-a]_q=\f{q^{a+b}-1-(q^{b-a}-1)}{q-1}=q^{b-a}\f{q^{2a}-1}{q-1}=q^{b-a}[2a]_q.$$ Then,
replacing the $(n+1)$-th column $C_{n+1}$ of the last determinant by $C_{n+1}-q^{b-a}\times C_{n}-[b-a]_q\times C_{1}$,
we see that
$$D_T(q)=\vmatrix {\bf 1}'&[d_q(v_j,v_k)]_{1\ls j,k\ls n-3}&[{\bs u}]_q'&[{\bs u}+a{\bf 1}]_q'&{\bf 0}'\\1&[{\bs u}]_q&0&[a]_q&0
\\1&[{\bs u}+a{\bf 1}]_q&[a]_q&0&f(a,b)
\\0&{\bf 0}&0&f(a,b)&g(a,b)
\\1-q&{\bf 1}&1&1&0\endvmatrix,$$
where
\begin{align*}g(a,b)=&\ (1-q^{b-a})f(a,b)-(1+q^{b-a})[a+b]_q
\\=&\ (q^{b-a}-q^{2(b-a)})[2a]_q-(1+q^{b-a})[a+b]_q.
\end{align*}
Now, expanding the last determinant according to its $n$th row, we find that
\begin{align*}D_T(q)=&\ f(a,b)\vmatrix {\bs 1}'&[d_q(v_j,v_k)]_{1\ls j,k\ls n-3}&[{\bs u}]_q'&{\bf 0}'\\1&[{\bs u}]_q&0&0
\\1&[{\bs u}+a{\bf 1}]_q&[a]_q&f(a,b)
\\1-q&{\bf 1}&1&0\endvmatrix
\\&\ -g(a,b)\vmatrix {\bf 1}'&[d_q(v_j,v_k)]_{1\ls j,k\ls n-3}&[{\bs u}]_q'&[{\bs u}+a{\bf 1}]_q'\\1&[{\bs u}]_q&0&[a]_q
\\1&[{\bs u}+a{\bf 1}]_q&[a]_q&0
\\1-q&{\bf 1}&1&1\endvmatrix
\end{align*}
and hence
\begin{align*}
D_T(q)=&\ -f(a,b)^2\vmatrix {\bs 1}'&[d_q(v_j,v_k)]_{1\ls j,k\ls n-3}&[{\bs u}]_q'\\1&[{\bs u}]_q&0
\\1-q&{\bf 1}&1\endvmatrix
\\&\ -g(a,b)\vmatrix {\bf 1}'&[d_q(v_j,v_k)]_{1\ls j,k\ls n-3}&[{\bs u}]_q'&[{\bs u}+a{\bf 1}]_q'\\1&[{\bs u}]_q&0&[a]_q
\\1&[{\bs u}+a{\bf 1}]_q&[a]_q&0
\\1-q&{\bf 1}&1&1\endvmatrix
\\=&\ -f(a,b)^2D_{T-\{v_{n-1},v_n\}}(q)-g(a,b)D_{T-\{v_n\}}(q)
\end{align*}
in view of \eqref{n-2} and \eqref{n-1}.
Combining this with the induction hypothesis for the trees $T-\{v_n,v_{n-1}\}$ and $T-\{v_n\}$, we obtain that
\begin{align*}D_T(q)=&\ -f(a,b)^2\prod_{e\in E(T)\sm\{e_1,e_2\}}[2w(e)]_q-g(a,b)\prod_{e\in E(T)\sm\{e_2\}}[2w(e)]_q
\\=&\ -q^{2(b-a)}[2a]_q^2\prod_{e\in E(T)\sm\{e_1,e_2\}}[2w(e)]_q-g(a,b)\prod_{e\in E(T)\sm\{e_2\}}[2w(e)]_q
\\=&\ (-q^{2b-2a}[2a]_q-g(a,b))\prod_{e\in E(T)\sm\{e_2\}}[2w(e)]_q.
\end{align*}
Observe that
\begin{align*}&-q^{2b-2a}[2a]_q-g(a,b)
\\=&\ -q^{2b-2a}[2a]_q-(q^{b-a}-q^{2(b-a)})[2a]_q+(1+q^{b-a})[a+b]_q
\\=&\ -\f{q^{b-a}(q^{2a}-1)}{q-1}+\f{(1+q^{b-a})(q^{a+b}-1)}{q-1}
\\=&\ \f{q^{b-a}-q^{a+b}}{q-1}+\f{q^{a+b}-1+q^{2b}-q^{b-a}}{q-1}=\f{q^{2b}-1}{q-1}=[2b]_q.
\end{align*}
Therefore
$$D_T(q)=[2b]_q\prod_{e\in E(T)\sm\{e_2\}}[2w(e)]_q=\prod_{e\in E(T)}[2w(e)]_q.$$

In view of the above, we have completed our induction proof of Theorem \ref{Th2.1}. \qed

\end{document}